\documentclass[12pt]{article}
\usepackage{latexsym,amsfonts,amssymb}
\usepackage[dvips]{color}
\usepackage{colordvi,multicol}

\def \cal{\mathcal}

\newtheorem{thm}{Theorem}[section]

\newtheorem{pro}[thm]{Proposition}
\newtheorem{defi}[thm]{Definition}
\newtheorem{rem}[thm]{Remark}
\newtheorem{exa}[thm]{Example}

\begin{document}
\title{\bf  A Note on Convergence of Random Variables}
\author{Ze-Chun Hu, Ting Ma\thanks{Corresponding
author: College of Mathematics, Sichuan University, Chengdu
610065,  China\vskip 0cm E-mail address: zchu@scu.edu.cn (Z.-C. Hu), matingting2008@scu.edu.cn (T. Ma), 1684680179@qq.com (X.-J. Zhu)}\,\, and Xiu-Ju Zhu\\ \\
 {\small College of Mathematics, Sichuan University}}

\maketitle
\date{}

\begin{abstract}
In this note,  convergence of random variables will be revisited. We will  give the answers to 5 questions among the  6 open questions introduced in (Convergence rates in the law of large numbers and new kinds of convergence of random variables,  {\it Communication in Statistics - Theory and Methods}, DOI: 10.1080/03610926.2020.1716248), and make some related discussions.
\end{abstract}



\noindent

\section{Introduction}\setcounter{equation}{0}

It is well known that all kinds of convergence of random variables play an important role in probability and statistics. In this note, we will make some discussions on convergence of random variables, and aim to give answers to 5 questions among the  6 open questions introduced in \cite{HS20}.

Let
$(\Omega,\mathcal{F},P)$ be a probability space and $\{X,X_n,n\geq
1\}$ be a sequence of random variables. We have the following kinds of convergence:

\begin{itemize}

\item $\{X_n,n\geq 1\}$  is said to almost surely converge to $X$,
if there exists a set $N\in \mathcal{F}$ such that $P(N)=0$ and
$\forall\omega\in \Omega\backslash N,\
\lim_{n\to\infty}X_n(\omega)=X(\omega)$, which is denoted by
$X_n\stackrel{a.s.}{\longrightarrow} X$ or $X_n\to X\ a.s.$.

\item $\{X_n,n\geq 1\}$  is said to converge to $X$ in probability, if for any $\varepsilon>0$, $\lim_{n\to\infty}P(\{|X_n-X|\geq \varepsilon\})=0$,  which is denoted by $X_n\stackrel{P}{\longrightarrow} X$.

\item $\{X_n,n\geq 1\}$  is said to $L^p$-converge to $X$ $(p>0)$
if $\lim_{n\to\infty}E[|X_n-X|^p]=0$, which is denoted by
$X_n\stackrel{L^p}{\longrightarrow} X$.

\item $\{X_n,n\geq 1\}$  is said to $L^{\infty}$-converge to $X$
if $\lim_{n\to\infty}\|X_n-X\|_{\infty}=0$, which is denoted by
$X_n\stackrel{L^{\infty}}{\longrightarrow} X$.

\item $\{X_n,n\geq 1\}$  is said to converge to $X$ in distribution, if for any bounded continuous function $f$, $\lim_{n\to\infty}E[f(X_n)]=E[f(X)]$,  which is denoted by $X_n\stackrel{d}{\longrightarrow} X$.

\item $\{X_n,n\geq 1\}$  is said to completely converge to $X$, if for any $\varepsilon>0$, $\sum_{n=1}^{\infty}P(\{|X_n-X|\geq \varepsilon\})<\infty$,  which is denoted by $X_n\stackrel{c.c.}{\longrightarrow} X$ (see \cite{HR47}).

\item $\{X_n,n\geq 1\}$  is said to S-$L^p$ converge to $X$
$(p>0)$ if $\sum_{n=1}^{\infty}E[|X_n-X|^p]<\infty$, which is
denoted by $X_n\stackrel{S\mbox{-}L^p}{\longrightarrow} X$ (see \cite[Definition 1.4]{LH17}).

\item $\{X_n,n\geq 1\}$  is said to strongly almost surely converge to
$X$ with order $\alpha$ ($\alpha>0$), if
$
\sum_{n=1}^{\infty}|X_n-X|^{\alpha}<\infty\ a.s.,
$
 which is denoted by $X_n\stackrel{S_{\alpha}\mbox{-}a.s.}{\longrightarrow} X$ (see
\cite[Definition 1.1]{HS20}).

\item $\{X_n,n\geq 1\}$  is said to strongly $L^{\infty}$-converge to
$X$ if
$
\sum_{n=1}^{\infty}\|X_n-X\|_{\infty}<\infty,
$
 which is denoted by $X_n\stackrel{S\mbox{-}L^{\infty}}{\longrightarrow} X$ (see
\cite[Definition 1.2]{HS20}).

\item $\{X_n,n\geq 1\}$  is said to $S_1\mbox{-}d$ converge to $X$, if
for any bounded Lipschitz continuous function $f$,
$
\sum_{n=1}^{\infty}|E[f(X_n)-f(X)]|<\infty,
$
which is denoted by $X_n\stackrel{S_1\mbox{-}d}{\longrightarrow}
X$ (see
\cite[Definition 1.3]{HS20}).

\item  $\{X_n,n\geq 1\}$  is said to $S_2\mbox{-}d$
converge to $X$, if for any continuous point $x$ of  $F$,
$
\sum_{n=1}^{\infty}|F_n(x)-F(x)|<\infty,
$
which is denoted by $X_n\stackrel{S_2\mbox{-}d}{\longrightarrow}
X$, where  $F_n$ and $F$ are the distribution functions of $X_n$ and $X$,
respectively (see
\cite[Definition 1.4]{HS20}).
\end{itemize}

In the final section of \cite{HS20}, the following 6 open questions were introduced:

\noindent {\bf Question 1.} What is the relation between the
$S_1\mbox{-}d$ convergence and the $S_2\mbox{-}d$ convergence?

\noindent {\bf Question 2.} Does
$X_n\stackrel{S\mbox{-}L^{\infty}}{\longrightarrow} X$ imply that
$X_n\stackrel{S_2\mbox{-}d}{\longrightarrow} X$?

\noindent {\bf Question 3.} Does
$X_n\stackrel{S\mbox{-}L^1}{\longrightarrow} X$ imply that
$X_n\stackrel{S_2\mbox{-}d}{\longrightarrow} X$?

\noindent {\bf Question 4.} Does
$X_n\stackrel{S_{\alpha}\mbox{-}a.s.}{\longrightarrow} X$
($\alpha>0$) imply that
$X_n\stackrel{S_1\mbox{-}d}{\longrightarrow} X$?

\noindent {\bf Question 5.} Does
$X_n\stackrel{c.c.}{\longrightarrow} X$ imply that
$X_n\stackrel{S_i\mbox{-}d}{\longrightarrow} X$ for $i\in \{1,2\}$?

\noindent {\bf Question 6.} Can we give a Skorokhod-type theorem
for the strong convergence in distribution and the
$S_{\alpha}\mbox{-}a.s.$ convergence?

In this note, we will give the answers to the first 5 questions, and make some related discussions. For simplicity, we introduce the following definition. 

\begin{defi}\label{def-2.4}
Let $\{X,X_n,n\geq 1\}$ be a sequence of random variables. If
for any bounded Lipschitz continuous function $f$, it holds that
$
\sum_{n=1}^{\infty}E[|f(X_n)-f(X)|]<\infty,
$
then $\{X_n,n\geq 1\}$  is said to $S_1^*\mbox{-}d$ converge to $X$, which is denoted by $X_n\stackrel{S_1^*\mbox{-}d}{\longrightarrow} X$.
\end{defi}

It is easy to know that  $X_n\stackrel{S\mbox{-}L^1}{\longrightarrow} X\Rightarrow X_n\stackrel{S_1^*\mbox{-}d}{\longrightarrow} X\Rightarrow X_n\stackrel{S_1\mbox{-}d}{\longrightarrow} X$.

Let $\{X,X_n,n\geq 1\}$ be a sequence of random variables. Denote by $\{f(t),f_n(t),n\geq 1\}$ the corresponding characteristic functions. It is well known that $X_n\stackrel{d}{\longrightarrow} X$ if and only if for any real number $t$,  $f_n(t)$ converges to $f(t)$ as $n\to\infty$. In virtue of this result, we introduce the following definition.

\begin{defi}\label{def-3.6}
Let $\{X,X_n,n\geq 1\}$ be a sequence of random variables. If
for any real number $t$, it holds that
$
\sum_{n=1}^{\infty}\left|E[e^{itX_n}]-E[e^{itX}]\right|<\infty,
$
then $\{X_n,n\geq 1\}$  is said to $S_3\mbox{-}d$ converge to $X$, which is denoted by $X_n\stackrel{S_3\mbox{-}d}{\longrightarrow} X$.
\end{defi}

It is easy to know that  $X_n\stackrel{S_1\mbox{-}d}{\longrightarrow} X\Rightarrow X_n\stackrel{S_3\mbox{-}d}{\longrightarrow} X$. Then we can rewrite the  diagram in \cite{HS20} as follows:
\begin{eqnarray*}
\begin{array}{ccccccc}
&& X_n\stackrel{S_1^*\mbox{-}d}{\longrightarrow} X & \Rightarrow  &  X_n\stackrel{S_1\mbox{-}d}{\longrightarrow} X & \Rightarrow & X_n\stackrel{S_3\mbox{-}d}{\longrightarrow} X\\
&&\Uparrow & &&&\Downarrow\\
X_n\stackrel{S\mbox{-}L^{\infty}}{\longrightarrow} X & \Rightarrow & X_n\stackrel{S\mbox{-}L^1}{\longrightarrow} X & \Rightarrow &
X_n\stackrel{S_1\mbox{-}a.s.}{\longrightarrow} X&&X_n\stackrel{d}{\longrightarrow} X\\
&&\Downarrow&&\Downarrow &&\Uparrow\\
&&X_n\stackrel{c.c.}{\longrightarrow} X& \Rightarrow& X_n\stackrel{a.s.}{\longrightarrow} X& \Rightarrow& X_n\stackrel{P}{\longrightarrow} X.\\
&&&&\Uparrow& & \Uparrow \\
&&&&X_n\stackrel{L^{\infty}}{\longrightarrow} X&\Rightarrow & X_n\stackrel{L^1}{\longrightarrow} X
\end{array}
\end{eqnarray*}

The rest of this note is organised as follows. In Section 2, we give the answers to Questions 1-5 based on 3 examples and some results in \cite{HS20}. In Section 3, we make more discussions on the relation among several kinds of convergence of random variables.

\section{Answers to Questions 1-5}\setcounter{equation}{0}

As to the answers to the first 5 questions introduced in Section 1, we have

\begin{pro}\label{pro-2.1}

(i) $X_n\stackrel{S_1\mbox{-}d}{\longrightarrow} X \nRightarrow X_n\stackrel{S_2\mbox{-}d}{\longrightarrow} X$, and $X_n\stackrel{S_2\mbox{-}d}{\longrightarrow} X \nRightarrow X_n\stackrel{S_1\mbox{-}d}{\longrightarrow} X$;

(ii) $X_n\stackrel{S\mbox{-}L^{\infty}}{\longrightarrow} X \nRightarrow X_n\stackrel{S_2\mbox{-}d}{\longrightarrow} X$;

(iii)$X_n\stackrel{S\mbox{-}L^1}{\longrightarrow} X \nRightarrow X_n\stackrel{S_2\mbox{-}d}{\longrightarrow} X$;

(iv) $X_n\stackrel{S_{\alpha}\mbox{-}a.s.}{\longrightarrow} X$
($\alpha>0$) $\nRightarrow$
$X_n\stackrel{S_1\mbox{-}d}{\longrightarrow} X$;

(v)$X_n\stackrel{c.c.}{\longrightarrow} X$ $\nRightarrow$
$X_n\stackrel{S_i\mbox{-}d}{\longrightarrow} X$ for $i\in \{1,2\}$.
\end{pro}

Before giving the proof of the above proposition, we present three examples.

\begin{exa}\label{exa-3.1}
 Let $\alpha>0$.
Define $\Omega=(0,1)$,  $\cal{F}=\cal{B}(\Omega)$ and $P$ be the
Lebesgue measure on $\Omega$. For $n\in\mathbb{N}$, we define a
random variable $X_n$ as follows:
\begin{eqnarray*}
X_n(\omega):=\left\{
\begin{array}{ll}
1, & \mbox{if}\ \   \omega\in (0,\frac{1}{n^2});\\
\frac{1}{n^{1/\alpha}}, & \mbox{if}\ \ \omega\in [\frac{1}{n^2},
1).
\end{array}
\right.
\end{eqnarray*}

By \cite[Example 3.11]{HS20}, we know that $X_n\stackrel{c.c.}{\longrightarrow} 0$, which together with \cite[Theorem 3.5]{HS20} implies that $X_n\stackrel{S_2\mbox{-}d}{\longrightarrow} 0$.

In the following, we will show that when $\alpha>1$,  $X_n\stackrel{S_1\mbox{-}d}{\nrightarrow} 0$. Let $f(x)=\sin x$. Then $f(x)$ is a bounded Lipschitz continuous function. We have
\begin{eqnarray*}
\sum_{n=1}^{\infty}|E[f(X_n)-f(0)]|&=&\sum_{n=1}^{\infty}|E[\sin X_n -\sin 0]|\nonumber\\
&=&\sum_{n=1}^{\infty}|E[\sin X_n]|\nonumber\\
&=&\sum_{n=1}^{\infty}\left[\frac{1}{n^2}\sin 1+\left(1-\frac{1}{n^2}\right)\sin\frac{1}{n^{1/\alpha}}\right]\nonumber\\
&=&\sin 1\sum_{n=1}^{\infty}\frac{1}{n^2}-\sum_{n=1}^{\infty}\frac{1}{n^2}\sin\frac{1}{n^{1/\alpha}}
+\sum_{n=1}^{\infty}\sin\frac{1}{n^{1/\alpha}}.
\end{eqnarray*}
It is easy to know that the first two sums are convergent.  By
$$
\lim_{n\to\infty}\frac{\sin\frac{1}{n^{1/\alpha}}}{\frac{1}{n^{1/\alpha}}}=1,
$$
and the fact that for $\alpha>1$,
$
\sum_{n=1}^{\infty}\frac{1}{n^{1/\alpha}}=\infty,
$
we know that the sum $\sum_{n=1}^{\infty}\sin\frac{1}{n^{1/\alpha}}$ is divergent. Hence
$$
\sum_{n=1}^{\infty}|E[f(X_n)-f(0)]|=\infty.
$$
It follows that $X_n\stackrel{S_1\mbox{-}d}{\nrightarrow} 0$.
\end{exa}

\begin{exa}\label{exa-3.2}
Define $\Omega=(0,1)$,  $\cal{F}=\cal{B}(\Omega)$ and $P$ be the
Lebesgue measure on $\Omega$. Let $\alpha,\beta$ be two constants satisfying $0<\alpha<1,\beta>1$. Let $X$ be a random variable defined on $(\Omega,\cal{F},P)$ with the density function $f(u)=(1-\alpha)(1-u)^{-\alpha},u\in (0,1)$. For any $n\in \mathbb{N}$, define a random variable
$$
X_n:=X+\frac{1}{n^{\beta}}.
$$
Then we have
$$
\sum_{n=1}^{\infty}\|X_n-X\|_{\infty}=\sum_{n=1}^{\infty}\frac{1}{n^{\beta}}<\infty,
$$
which implies that $X_n\stackrel{S\mbox{-}L^{\infty}}{\longrightarrow} X$.

 Denote by $F_n$ and $F$ the distribution functions of $X_n$ and $X$, respectively. Suppose that $(1-\alpha)\beta\leq 1$. Then we have
\begin{eqnarray*}
\sum_{n=1}^{\infty}|F_n(1)-F(1)|&=&\sum_{n=1}^{\infty}\left|F\left(1-\frac{1}{n^{\beta}}\right)-F(1)\right|\\
&=&\sum_{n=1}^{\infty}\left[F(1)-F\left(1-\frac{1}{n^{\beta}}\right)\right]\\
&=&\sum_{n=1}^{\infty}\int_{1-\frac{1}{n^{\beta}}}^1(1-\alpha)(1-u)^{-\alpha}du\\
&=&(1-\alpha)\sum_{n=1}^{\infty}\int_{0}^{\frac{1}{n^{\beta}}}v^{-\alpha}dv\\
&=&\sum_{n=1}^{\infty}\frac{1}{n^{(1-\alpha)\beta}}=\infty.
\end{eqnarray*}
It follows that $X_n\stackrel{S_2\mbox{-}d}{\nrightarrow} X$.
\end{exa}

\begin{exa}\label{exa-3.3}
Define $\Omega=(0,1)$,  $\cal{F}=\cal{B}(\Omega)$ and $P$ be the
Lebesgue measure on $\Omega$. For $n\in\mathbb{N}$, we define a
random variable $X_n$ as follows:
\begin{eqnarray*}
X_n(\omega):=\left\{
\begin{array}{ll}
1, & \mbox{if}\ \   \omega\in (0,\frac{1}{n});\\
0, & \mbox{if}\ \ \omega\in [\frac{1}{n},1).
\end{array}
\right.
\end{eqnarray*}
It is easy to check that for any $\alpha>0$, any $\omega\in (0,1)$, we have
$$
\sum_{n=1}^{\infty}|X_n(\omega)-0|^{\alpha}<\infty.
$$
It follows that $X_n\stackrel{S_{\alpha}\mbox{-}a.s.}{\rightarrow} 0$.

Let $f(x)=\sin x$. Then $f$ is a bounded Lipschitz continuous function. We have
\begin{eqnarray*}
\sum_{n=1}^{\infty}|E[f(X_n)-f(0)]|&=&\sum_{n=1}^{\infty}|E[\sin X_n-\sin 0]|\\
&=&\sum_{n=1}^{\infty}|E[\sin X_n]|\\
&=&\sum_{n=1}^{\infty}\left|\frac{1}{n}\sin 1+\left(1-\frac{1}{n}\right)\sin 0\right|\\
&=&\sin 1\sum_{n=1}^{\infty}\frac{1}{n}\\
&=&\infty,
\end{eqnarray*}
and thus $X_n\stackrel{S_1\mbox{-}d}{\nrightarrow} 0$.
\end{exa}

\noindent {\bf Proof of Proposition \ref{pro-2.1}: }

 (i) By Example \ref{exa-3.1}, we get that $X_n\stackrel{S_2\mbox{-}d}{\longrightarrow} X\nRightarrow X_n\stackrel{S_1\mbox{-}d}{\longrightarrow} X$. By Example \ref{exa-3.2} and the fact that $X_n\stackrel{S\mbox{-}L^{\infty}}{\longrightarrow} X\Rightarrow X_n\stackrel{S\mbox{-}L^1}{\longrightarrow} X\Rightarrow X_n\stackrel{S_1\mbox{-}d}{\longrightarrow} X$, we obtain that $X_n\stackrel{S_1\mbox{-}d}{\longrightarrow} X \nRightarrow X_n\stackrel{S_2\mbox{-}d}{\longrightarrow} X$.

(ii) By Example \ref{exa-3.2}, we obtain that $X_n\stackrel{S\mbox{-}L^{\infty}}{\longrightarrow} X \nRightarrow X_n\stackrel{S_2\mbox{-}d}{\longrightarrow} X$.

(iii) By Example \ref{exa-3.2} and the fact that $X_n\stackrel{S\mbox{-}L^{\infty}}{\longrightarrow} X\Rightarrow X_n\stackrel{S\mbox{-}L^1}{\longrightarrow} X$, we obtain that $X_n\stackrel{S\mbox{-}L^1}{\longrightarrow} X \nRightarrow X_n\stackrel{S_2\mbox{-}d}{\longrightarrow} X$.

(iv) By Example \ref{exa-3.3}, we obtain that $X_n\stackrel{S_{\alpha}\mbox{-}a.s.}{\longrightarrow} X$
($\alpha>0$) $\nRightarrow$
$X_n\stackrel{S_1\mbox{-}d}{\longrightarrow} X$.

(v) By Example \ref{exa-3.1}, we get that $X_n\stackrel{c.c.}{\longrightarrow} X \nRightarrow X_n\stackrel{S_1\mbox{-}d}{\longrightarrow} X$. 
 By Example \ref{exa-3.2} and the fact that $X_n\stackrel{S\mbox{-}L^{\infty}}{\longrightarrow} X\Rightarrow X_n\stackrel{S\mbox{-}L^1}{\longrightarrow} X\Rightarrow X_n\stackrel{c.c.}{\longrightarrow} X$, we obtain that $X_n\stackrel{c.c.}{\longrightarrow} X \nRightarrow X_n\stackrel{S_2\mbox{-}d}{\longrightarrow} X$. \hfill\fbox

\section{More discussions}

In this section, we make more discussions on the relations of several kinds of convergence of random variables.

\subsection{Main results and questions}

By Proposition \ref{pro-2.1}(ii),  we know that generally speaking, $X_n\stackrel{S\mbox{-}L^{\infty}}{\longrightarrow} X$ does not imply $X_n\stackrel{S_2\mbox{-}d}{\longrightarrow} X$. But if some addition condition is assumed, we may have that $X_n\stackrel{S\mbox{-}L^{\infty}}{\longrightarrow} X\Rightarrow X_n\stackrel{S_2\mbox{-}d}{\longrightarrow} X$. By \cite[Proposition 3.7(i)]{HS20} and the fact that $X_n\stackrel{S\mbox{-}L^{\infty}}{\longrightarrow} X\Rightarrow X_n\stackrel{S\mbox{-}L^1}{\longrightarrow} X\Rightarrow X_n\stackrel{c.c.}{\longrightarrow} X$, we get that if $X$ is a discrete random variable such that $\{x\in\mathbb{R}:P(X=x)=0\}$ is an open subset of
$\mathbb{R}$ and $X_n\stackrel{S\mbox{-}L^{\infty}}{\longrightarrow} X$, then $X_n\stackrel{S_2\mbox{-}d}{\longrightarrow} X$. The following proposition extends this result.

\begin{pro}\label{pro-2.2}
Let $\{X,X_n,n\geq 1\}$ be a sequence of random variables and
$\{F,F_n,n\geq 1\}$ be the corresponding sequence of distribution
functions.  If $F$ is locally Lipschitz continuous at each continuous point $x$ of $F$ and $
X_n\stackrel{S\mbox{-}L^{\infty}}{\longrightarrow} X$, then $X_n\stackrel{S_2\mbox{-}d}{\longrightarrow} X.$
\end{pro}

\begin{rem}\label{lem-2.3}
(i) By \cite[Theorem 3.5]{HS20}, we know that if $C$ is a constant, then $X_n\stackrel{c.c.}{\longrightarrow} C\Leftrightarrow X_n\stackrel{S_2\mbox{-}d}{\longrightarrow} C.$

(ii) By \cite[Proposition 3.7]{HS20}, we know that if $X_n\stackrel{c.c.}{\longrightarrow} X$, and $X$ is a discrete random variable such that $\{x\in \mathbb{R}: P(X=x)=0\}$ is an open subset of $\mathbb{R}$, then $X_n\stackrel{S_2\mbox{-}d}{\longrightarrow} X$.

(iii) By (ii) and the fact that $X_n\stackrel{S\mbox{-}L^1}{\longrightarrow} X\Rightarrow X_n\stackrel{c.c.}{\longrightarrow} X$, we get that if $X_n\stackrel{S\mbox{-}L^1}{\longrightarrow} X$, and $X$ is a discrete random variable such that $\{x\in \mathbb{R}: P(X=x)=0\}$ is an open subset of $\mathbb{R}$, then $X_n\stackrel{S_2\mbox{-}d}{\longrightarrow} X$.
\end{rem}

\begin{pro}\label{pro-2.5}
If $X_n\stackrel{c.c.}{\longrightarrow} X$ and  $\sum_{n=1}^{\infty}E[|X_n-X|I_{\{|X_n-X|<\varepsilon\}}]<\infty$ for some positive number $\varepsilon$, then $X_n\stackrel{S_1^*\mbox{-}d}{\longrightarrow} X$, and thus $X_n\stackrel{S_1\mbox{-}d}{\longrightarrow} X$.
\end{pro}

\begin{rem}\label{rem-2.6}
The condition that $\sum_{n=1}^{\infty}E[|X_n-X|I_{\{|X_n-X|<\varepsilon\}}]<\infty$ for some positive number $\varepsilon$ is necessary in the sense that if $X_n\stackrel{c.c.}{\longrightarrow} X$ and $X$ is a bounded random variable, then $X_n\stackrel{S_1^*\mbox{-}d}{\longrightarrow} X$ if and only if $\sum_{n=1}^{\infty}E[|X_n-X|I_{\{|X_n-X|<\varepsilon\}}]<\infty$ for any positive number $\varepsilon$.

Obviously, by virtue of Proposition \ref{pro-2.5}, we need only to show the necessity. Suppose that $X_n\stackrel{c.c.}{\longrightarrow} X$, $X$ is a bounded random variable, and $X_n\stackrel{S_1^*\mbox{-}d}{\longrightarrow} X$. Let $M$ be a positive number satisfying $|X|\leq M$ a.s.. For any positive number $\varepsilon$, define a function as follows:
\begin{eqnarray*}
f_\varepsilon(x):=\left\{
\begin{array}{cl}
M+\varepsilon& \mbox{if}\ x>M+\varepsilon;\\
x& \mbox{if}\ -M-\varepsilon\leq x\leq M+\varepsilon;\\
-M-\varepsilon& \mbox{if}\ x<-M-\varepsilon.
\end{array}
\right.
\end{eqnarray*}
It is easy to check that $f_\varepsilon$ is a bounded Lipschitz continuous function. By the definitions of $X_n\stackrel{S_1^*\mbox{-}d}{\longrightarrow} X$ and $f_\varepsilon$, we have
\begin{eqnarray*}
\infty&>&\sum_{n=1}^{\infty}E[|f_\varepsilon(X_n)-f_\varepsilon(X)|]\\
&=&\sum_{n=1}^{\infty}E[|f_\varepsilon(X_n)-f_\varepsilon(X)|I_{\{|X_n-X|<\varepsilon\}}]
+\sum_{n=1}^{\infty}E[|f_\varepsilon(X_n)-f_\varepsilon(X)|I_{\{|X_n-X|\geq \varepsilon\}}]\\
&\geq&\sum_{n=1}^{\infty}E[|f_\varepsilon(X_n)-f_\varepsilon(X)|I_{\{|X_n-X|<\varepsilon\}}]\\
&=&\sum_{n=1}^{\infty}E[|X_n-X|I_{\{|X_n-X|<\varepsilon\}}].
\end{eqnarray*}
\end{rem}

In virtue of the $S_3\mbox{-}d$ convergence, we have the following three questions:

\noindent {\bf Question 7.} What is the relation between the
$S_2\mbox{-}d$ convergence and the $S_3\mbox{-}d$ convergence? 

\noindent {\bf Question 8.} Does
$X_n\stackrel{S_{\alpha}\mbox{-}a.s.}{\longrightarrow} X$
($\alpha>0$) imply that
$X_n\stackrel{S_3\mbox{-}d}{\longrightarrow} X$?

\noindent {\bf Question 9.} Does
$X_n\stackrel{c.c.}{\longrightarrow} X$ imply that
$X_n\stackrel{S_3\mbox{-}d}{\longrightarrow} X$?

\bigskip 

Since $X_n\stackrel{S_1\mbox{-}d}{\longrightarrow} X \nRightarrow X_n\stackrel{S_2\mbox{-}d}{\longrightarrow} X$ and $X_n\stackrel{S_1\mbox{-}d}{\longrightarrow} X \Rightarrow X_n\stackrel{S_3\mbox{-}d}{\longrightarrow} X$, we get that $X_n\stackrel{S_3\mbox{-}d}{\longrightarrow} X \nRightarrow X_n\stackrel{S_2\mbox{-}d}{\longrightarrow} X$.   By checking Example \ref{exa-3.1}, we can get that $X_n\stackrel{S_2\mbox{-}d}{\longrightarrow} X \nRightarrow X_n\stackrel{S_3\mbox{-}d}{\longrightarrow} X$.

By checking Example \ref{exa-3.3}, we get that $X_n\stackrel{S_{\alpha}\mbox{-}a.s.}{\longrightarrow} X$
($\alpha>0$) $\nRightarrow$
$X_n\stackrel{S_3\mbox{-}d}{\longrightarrow} X$.

By checking Example \ref{exa-3.1}, we get that $X_n\stackrel{c.c.}{\longrightarrow} X \nRightarrow X_n\stackrel{S_3\mbox{-}d}{\longrightarrow} X$.

\begin{rem}
\cite[Example 3.12]{HS20} shows that if X is nondegenerate, then $X_n\stackrel{S_3\mbox{-}d}{\longrightarrow} X$ does not imply $X_n\stackrel{P}{\longrightarrow} X$.
\end{rem}

\subsection{Proofs}\setcounter{equation}{0}

\noindent {\bf Proof of Proposition \ref{pro-2.2}.} Suppose that $X_n\stackrel{S\mbox{-}L^{\infty}}{\longrightarrow} X$. Denote $\alpha_n=\|X_n-X\|_{\infty}$. Then $\alpha_n\geq 0$ and
\begin{eqnarray}\label{pro-2.5-a}
\sum_{n=1}^{\infty}\alpha_n<\infty.
\end{eqnarray}
For any $x\in \mathbb{R}$, we have
\begin{eqnarray*}
F_n(x)-F(x)&=&P(X_n\leq x)-F(x)\\
&=&P(X+X_n-X\leq x)-F(x)\\
&\leq& P(X\leq x+\alpha_n)-F(x)\\
&=&F(x+\alpha_n)-F(x),
\end{eqnarray*}
and
\begin{eqnarray*}
F_n(x)-F(x)&=&1-P(X_n> x)-F(x)\\
&=&1-P(X+X_n-X> x)-F(x)\\
&\geq&1-P(X>x-\alpha_n)-F(x)\\
&=&F(x-\alpha_n)-F(x)\\
&=&-[F(x)-F(x-\alpha_n)].
\end{eqnarray*}
It follows that
\begin{eqnarray}\label{pro-2.5-b}
|F_n(x)-F(x)|\leq (F(x+\alpha_n)-F(x))+(F(x)-F(x-\alpha_n)).
\end{eqnarray}

Let $x$ be any continuous point of $F$. By the assumption, there exist two constants $K$ and $\delta$ such that for any $u,v\in (x-\delta,x+\delta)$,
\begin{eqnarray}\label{pro-2.5-b-1}
|F(u)-F(v)|\leq K|u-v|.
\end{eqnarray}
Since $\lim_{n\to\infty}\alpha_n=0$, there exists $N_0$ such that for any $n> N_0$, we have $\alpha_n<\delta$. Then by (\ref{pro-2.5-b-1}), we get that for any $n> N_0$,
$$
F(x+\alpha_n)-F(x)\leq K\alpha_n.
$$
Then by (\ref{pro-2.5-a}), we get that
\begin{eqnarray}\label{pro-2.5-c}
\sum_{n=1}^{\infty}(F(x+\alpha_n)-F(x))\leq \sum_{n=1}^{N_0}(F(x+\alpha_n)-F(x))+K\sum_{n=N_0+1}^{\infty}\alpha_n<\infty.
\end{eqnarray}
Similarly, we have
\begin{eqnarray}\label{pro-2.5-d}
\sum_{n=1}^{\infty}(F(x)-F(x-\alpha_n))<\infty.
\end{eqnarray}
By (\ref{pro-2.5-b}), (\ref{pro-2.5-c}) and (\ref{pro-2.5-d}), we get that
$$
\sum_{n=1}^{\infty}|F_n(x)-F(x)|<\infty.
$$
It follows that $X_n\stackrel{S_2\mbox{-}d}{\longrightarrow} X.$ \hfill\fbox

\bigskip

\noindent {\bf Proof of Proposition \ref{pro-2.5}.} Suppose that $f$ is a bounded Lipschitz continuous function. Then there exist two constants $K,M$ such that for any $x,y\in \mathbb{R}$, we have $|f(x)|\leq M$ and $|f(x)-f(y)|\leq K|x-y|$. Then by the assumptions, we get
\begin{eqnarray*}
\sum_{n=1}^{\infty}E[|f(X_n)-f(X)|]&=&\sum_{n=1}^{\infty}E[|f(X_n)-f(X)|I_{\{|X_n-X|<\varepsilon\}}]+
\sum_{n=1}^{\infty}E[|f(X_n)-f(X)|I_{\{|X_n-X|\geq \varepsilon\}}]\\
&\leq& K\sum_{n=1}^{\infty}E[|X_n-X|I_{\{|X_n-X|<\varepsilon\}}]+
2M\sum_{n=1}^{\infty}P(|X_n-X|\geq \varepsilon)\\
&<&\infty.
\end{eqnarray*}
Thus $X_n\stackrel{S_1\mbox{-}d}{\longrightarrow} X$.\hfill\fbox
\bigskip

\bigskip

{ \noindent {\bf\large Acknowledgments} \quad  We are grateful to the support of NNSFC (Grant No. 11771309 and No. 11871184).}

\end{document}